# A Digital Twin Based Decision Support System for the Management of an Operating Room

Leah RIFI
Industrial Engineering Center, IMT Mines Albi, Albi, Tarn, France, email: leah.rifi@mines-albi.fr

Canan PEHLIVAN
Industrial Engineering Center, IMT Mines Albi, Albi, Tarn, France, email: canan.pehlivan@mines-albi.fr

Cléa MARTINEZ
Industrial Engineering Center, IMT Mines Albi, Albi, Tarn, France, email: clea.martinez@mines-albi.fr

Maria DI MASCOLO
Univ. Grenoble Alpes, CNRS, Grenoble INP, G-SCOP, 38000 Grenoble, France, email: maria.di-mascolo@grenoble-inp.fr

Franck FONTANILI
Industrial Engineering Center, IMT Mines Albi, Albi, Tarn, France, franck.fontanili@mines-albi.fr

With healthcare demand rising worldwide, hospital services are increasingly needed. Hospitals' performance is tightly linked to their surgical suite performance, which makes it necessary for surgical suites to be efficient. In this paper, we focus on the operating room schedule execution and related decision-making. We thus propose a digital twin-based decision support system for the prospective and retrospective simulation and analysis of the operating room schedule execution. We describe the developed prototype (inputs/output/parameters/modeling) and its functionalities and show its application to an operating room inspired by a real case.

## 1 Context and problematic

As healthcare demand escalates globally, the necessity for medical services grows, and the efficiency of hospitals' surgical suites becomes pivotal. These suites not only represent a significant portion of hospital budgets, with 40%-50% allocated to them, but they also constitute around 30% of overall healthcare costs (Kaye et al., 2020; Macario et al., 1997). Thus, ensuring the effectiveness and profitability of operating rooms (OR) is imperative. However, the management of a surgical suite is complex due to several factors: the **variety of patient pathways** (outpatient and inpatient, elective and non-elective patients), the **diversity of professionals** working within the suite (surgeons, anesthesiologists, nurses, etc.), the **interactions with internal and external units** (wards, diagnostic facilities, external service providers…), and the **uncertainties inherent in medical practice** (variability in procedure durations and unexpected events such as cancellations and emergency arrivals…) (Zhu et al., 2019).

Against this backdrop, the objective of this research is to advocate for a performant, resilient, robust organization of surgical suites. Performance is evaluated through key metrics such as staff workload (overtime), patient satisfaction (waiting time), and operational efficiency (operation room utilization). A robust organization can maintain performance levels despite inherent variability, while a resilient one can withstand significant disruptions without compromising performance. Achieving this objective not only enhances efficiency and profitability but also contributes to delivering safe and high-quality medical care.

To attain this objective, the study focuses on optimizing the operating schedule, known as the OR planning and scheduling problem. This involves decision-making across four hierarchical levels: strategic, tactical, offline operational, and online operational. While existing literature primarily addresses strategic, tactical, and offline operational planning, there's a notable gap in addressing disruptions at the online operational level, crucial for real-time decision-making (Guerriero & Guido, 2011).

Within this context, we investigate the following general research question: "How can we promote and maintain the performance of a surgical suite's organization under uncertainties?". To answer this question, we propose: (1) A methodology for a prospective analysis of the provisional schedule, and (2) A methodology for a retrospective analysis of the performed schedule. Leveraging insights from a digital twin-based decision support system (DT-DSS), the study seeks to improve decision-making processes for OR managers, facilitating both proactive schedule preparation and post-operative analysis. A case study based on the Hôpital Privé de La Baie (GIE Vivalto Santé, France) serves as a proof of concept, demonstrating the applicability of the proposed methodologies and their potential to enhance operational efficiency across various surgical suites.

The remainder of this article is divided into two parts: our solution proposal and an illustrative case.

## 2  Solution Proposal

### 2.1  Prospective and Retrospective Analyses Description

We propose a prospective approach to analyzing the provisional schedule, that is divided into five key steps. The objectives are to (1) anticipate and address potential disruptions before the actual execution of the schedule, (2) to evaluate the effectiveness of predictive disruption management, and (3) improve the surgery scheduling and sequencing. First, we assess the practicality and feasibility of the provisional schedule (step #1). Then, we evaluate the provisional schedule's performance (step #2). Following this, we analyze the impact of stochastic durations on schedule performance to gauge its robustness against uncertainties (step #3). Subsequently, we introduce additional non-elective arrivals and examine the schedule's performance under different scheduling strategies to assess the resilience of the organization (step #4). Finally, we simulate the execution of the provisional schedule within a stochastic environment incorporating both variable durations and non-elective arrivals to assess the performance, the robustness and the resilience of the provisional schedule simultaneously (step #5).

On the other hand, we propose a retrospective analysis of the performed schedule. Here, the aim is to evaluate the quality of reactive disruption management and to identify whether performance shortcomings originate from offline or online decisions. To begin with, we assess the respect of resource constraints in a deterministic environment (step #1). Then, we evaluate the performed schedule's performance (step #2). Finally, we analyze any disparities between the performed and provisional schedules to determine whether they stem from offline or online operational decisions (step #3).

### 2.2  Prototype Description

In this section, we present the DT-DSS developed prototype, outlining its inputs, outputs, parameters, and modeling approach. The input data comprises a relational database consisting of three tables: the master surgery schedule, the attributes of scheduled cases, and the durations of scheduled cases. The output is represented through a dashboard, showcasing key performance indicators such as staff overtime, patient waiting time, and operating room utilization. Additionally, the dashboard includes a Gantt Chart Diagram illustrating the status of each OR, including idle time, resource allocation, setup, and procedure duration.

The parameters of the prototype encompass various aspects, including the initial schedule type, the process type, the constraints on resources, the duration type, whether to keep or not the non-elective cases of the initial schedule, whether to add non-elective arrivals to the initial schedule, and the number of replications.

The Digital Twin-based Decision Support System (DT-DSS) is constructed utilizing a modeling and simulation tool, Flexsim Healthcare®, which facilitates the modeling of material resources, human resources, and processes within the surgical suite. This tool enables the simulation of schedule execution, providing insights into the system's performance under different scenarios.

## 3  Application to an Operating Room Inspired by A Real-Case

In this section, we present a synthesis of both Step #1 and Step #2 of the retrospective analysis, focusing on a one-day use case inspired by our hospital partner.

In Figure 1, we display the OR state Gantt Chart Diagram of the performed schedule for different scheduling strategies (real-life, first fit, worst fit, and best fit). Each row corresponds to an OR and each color corresponds to an OR state: setup with anesthesiologist (yellow), setup without anesthesiologist (blue), surgical procedure (green), reversal (purple), idle time (light gray), and off schedule (dark grey). We indicate non-elective cases with a red or an orange line. For better readability, we display all the ORs for the real-life scenario, but only the ORs with non-elective cases for the other scenarios.

During the real performed schedule execution (step #1), the OR manager dealt with two urgent cases and scheduled them in OR#1 and #2. Both cases were entered at the earliest time possible and were the last case in their respective OR. In step #2, we simulate the performed schedule execution as it happened in real-life, but we change the urgent case scheduling solution. BF and WF strategies give the same scheduling solution because there are no shift end durations available for any of the urgent cases. We note that the FF strategy leads to scheduling the case in an OR#2 that has been closed for several hours, and that the BF strategy leads to scheduling both cases one after the other.

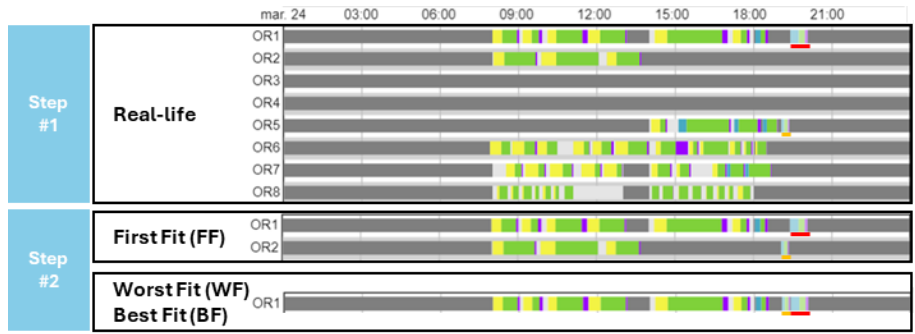

**Figure 1 – Real-life scheduling of the non-elective (step #1) VS other scheduling solutions (step #2)**

In Table 1, we display the value of two of our KPI for each scenario. The OR utilization and the staff overtime are the same for all scenarios because the urgent cases had arrival times and preoperative durations such that they could only be scheduled in overtime, after the last OR shift end. We note that both KPI do not reach their performance target. The prospective analysis results also offer insight into whether the prospective analysis was relevant or not.

**Table 1 - KPI value for each scenario**

| KPI | Real-life | FF | BF / WF | Rate Target |
|---|---|---|---|---|
| OR Utilization (percentage of OR shift) | 77.3% | 77.3% | 77.3% | ~85% |
| Staff Overtime (percentage of OR shift) | 9.6% | 9.6% | 9.6% | ≤ 5% |

The information gathered during these two steps can help support organizational decisions such as opening or closing OR (punctually or regularly) or serve to either encourage or discourage OR managers to implement specific non-elective scheduling strategies.